\RequirePackage{ifpdf}
\ifpdf 
\documentclass[pdftex]{sigma}
\else
\documentclass{sigma}
\fi

\begin{document}

\allowdisplaybreaks

\renewcommand{\thefootnote}{$\star$}

\renewcommand{\PaperNumber}{042}

\FirstPageHeading

\ShortArticleName{A Lax Formalism for the Elliptic Dif\/ference Painlev\'e Equation}

\ArticleName{A Lax Formalism for the Elliptic Dif\/ference\\ Painlev\'e Equation\footnote{This paper is a contribution to the Proceedings of the Workshop ``Elliptic Integrable Systems, Isomonodromy Problems, and Hypergeometric Functions'' (July 21--25, 2008, MPIM, Bonn, Germany). The full collection
is available at
\href{http://www.emis.de/journals/SIGMA/Elliptic-Integrable-Systems.html}{http://www.emis.de/journals/SIGMA/Elliptic-Integrable-Systems.html}}}

\Author{Yasuhiko YAMADA}

\AuthorNameForHeading{Y. Yamada}

\Address{Department of Mathematics, Faculty of Science,
Kobe University, Hyogo 657-8501, Japan}

\Email{\href{mailto:yamaday@math.kobe-u.ac.jp}{yamaday@math.kobe-u.ac.jp}}

\ArticleDates{Received November 20, 2008, in f\/inal form March 25,
2009; Published online April 08, 2009}

\Abstract{A Lax formalism for the elliptic Painlev\'e equation
is presented. The construction is based on the geometry of the
curves on ${\mathbb P}^1\times{\mathbb P}^1$ and described
in terms of the point conf\/igurations.}


\Classification{34A05; 14E07; 14H52}

\renewcommand{\thefootnote}{\arabic{footnote}}
\setcounter{footnote}{0}

\section{Introduction}

The dif\/ference analogs of the Painlev\'e dif\/ferential equations
have been extensively studied in the last two
decades (see \cite{GNR,Sakai} for example).
It is now widely recognized that some of the aspects of the
Painlev\'e equations, in particular their algebraic or geometric
properties, can be understood in universal
way by considering dif\/ferential and dif\/ference cases together.

In \cite{Sakai}, Sakai studied the dif\/ference Painlev\'e equations
from the point of view of rational surfaces and classif\/ied them
into three categories:
additive, multiplicative ($q$-dif\/ference) and elliptic\footnote{The addition formulae of the trigonometric/elliptic functions are
the typical examples of the multiplicative/elliptic dif\/ference
equations.}.
The classif\/ication is summarized in the following diagram:
{\arraycolsep=1pt
\[
\begin{array}{ccccccccccccccccccc}
{\rm ell.}\quad&E^{(1)}_8\\[-2mm]
&&&&&&&&&&&&&&&{\mathbb Z}\\
&&&&&&&&&&&&&&\nearrow\\
{\rm mul.}\quad&E^{(1)}_8&\rightarrow&E^{(1)}_7&\rightarrow&E^{(1)}_6
&\rightarrow &D^{(1)}_5&\rightarrow&A^{(1)}_4
&\rightarrow&(A_2\!+\!A_1)^{(1)}&\rightarrow&(A_1\!+\!A_1)^{(1)}
&\rightarrow&A^{(1)}_1&\rightarrow&{\mathcal D}_6\\[2mm]
\\
{\rm add.}\quad&E^{(1)}_8&\rightarrow&E^{(1)}_7&\rightarrow&E^{(1)}_6
&&\rightarrow&&D^{(1)}_4&\rightarrow&A_3^{(1)}&\rightarrow
&(A_1\!+\!A_1)^{(1)}&\rightarrow&A^{(1)}_1&\rightarrow&{{\mathbb Z}_2}\\
&&&&&&&&&&&&\searrow&&\searrow&&&\downarrow\\
&&&&&&&&&&&&&A^{(1)}_2&\rightarrow&A^{(1)}_1&\rightarrow&1\\
\end{array}
\]

}

Among them, Sakai's elliptic Painlev\'e equation \cite{Sakai}
is the master equation of all the second order Painlev\'e equations.
It has the af\/f\/ine Weyl group symmetry of type $E^{(1)}_8$ and all
the other cases arise as its degenerations.

It is well known that the dif\/ferential Painlev\'e equations describe
iso-monodromy deformations of linear dif\/ferential equations.
Since the iso-monodromy interpretation of the Painlev\'e equations
is a main source of variety of deep properties of the latter, it is
an important problem to f\/ind Lax formalisms for dif\/ference cases.

In fact, for some of dif\/ference Painlev\'e equations, the Lax
formulations have been known
(see \cite{AB,AL,B1,B,JS} for example).
Let us give an example of $q$-dif\/ference case with symmetry of
type~$D^{(1)}_5$. The equation is Jimbo--Sakai's $q$-$P_{\rm VI}$
equation \cite{JS} which is a discrete dynamical system def\/ined
by the following birational transformation:
\begin{gather}
T: \ \left(
\begin{array}{cccc}
a_1,&a_2,&a_3,&a_4\\
b_1,&b_2,&b_3,&b_4
\end{array};f,g \right)
\mapsto
\left(
\begin{array}{cccc}
qa_1,&qa_2,&a_3,&a_4\\
qb_1,&qb_2,&b_3,&b_4
\end{array}; \dot{f}, \dot{g} \right),\nonumber\\
\dot{f}{f}
=\frac{(\dot{g}-b_1)(\dot{g}-b_2)}
{(\dot{g}-b_3)(\dot{g}-b_4)}a_3 a_4, \qquad
\dot{g}{g}
=\frac{({f}-a_1)({f}-a_2)}
{({f}-a_3)({f}-a_4)}b_3b_4,\label{eq:qP6}
\end{gather}
where $(f,g) \in {\mathbb P}^1\times{\mathbb P}^1$ are the dependent
variables and
$a_1,\ldots,a_4$, $b_1,\ldots,b_4$ are complex parame\-ters with a constraint
$q=a_3a_4 b_1 b_2/(a_1 a_2 b_3 b_4)$.

The $q$-$P_{\rm VI}$ equation \eqref{eq:qP6} was originally derived as
the compatibility of certain $2\times 2$ matrix Lax pair:
\begin{gather}
Y(qz)=A(z) Y(z), \qquad
T(Y(z))=B(z) Y(z),\nonumber
\end{gather}
which is equivalent (up to a gauge transformation) with the
following scalar Lax pair for the f\/irst component $y$ of $Y$.
One of the scalar Lax equations is
\begin{gather}
\frac{(a_1-z)(a_2-z)}{a_1a_2(z-f)}y(qz)-\left(
c_0+c_1 z+\frac{c_2 z}{z-{f}}+\frac{c_3 z}{z-q {f}}\right)y(z)
\nonumber\\
\qquad{}
+\frac{a_1a_2(z-qa_3)(z-qa_4)}{b_3b_4 q^2(z-qf)}y\left(\frac{z}{q}\right)=0
\label{eq:qP6L1},
\end{gather}
where
$c_0=-\frac{a_1 a_2}{f}\big(\frac{1}{b_1}+\frac{1}{b_2}\big)$,
$c_1=\frac{1}{q}\big(\frac{1}{b_3}+\frac{1}{b_4}\big)$,
$c_2=\frac{(f-a_1)(f-a_2)}{qfg}$ and
$c_3=\frac{(f-a_3)(f-a_4)g}{b_3 b_4 f}$. The other one is
\begin{gather}
q g y(qz)-a_1a_2y(z)+z(z-f)T^{-1}(y(z))=0\label{eq:qP6L2}.
\end{gather}

In elliptic case, it is natural to expect a scalar Lax pair which looks like
\begin{gather}
C_1 y(z-\delta)+C_2 y(z)+C_3 y(z+\delta)=0,\qquad
C_4 y(z-\delta)+C_5 y(z)+C_6 T(y(z))=0,\label{eq:L12gen}
\end{gather}
where the coef\/f\/icients $C_1,\ldots,C_6$ are elliptic
functions on variables $z$ and other parameters. However,
the explicit construction of such
Lax formalism has remained as a dif\/f\/icult problem
because of the complicated elliptic dependence including many parameters.

In this paper, we give a Lax formulation for the elliptic Painlev\'e
equation with $E^{(1)}_8$ symmetry using a geometric method.
The main idea is to consider the Lax pair~\eqref{eq:L12gen} as
equations for algebraic curves with respect to the unknown
variables of the Painlev\'e equation.
We note that the Lax pair for the additive dif\/ference $E^{(1)}_8$ case
was obtained by Boalch~\cite{B1}. For the elliptic case, another approach
has recently proposed by Arinkin, Borodin and Rains \cite{ABR,R}.

This paper is organised as follows.
In Section~\ref{sect:ell-p}, a geometric description of the
elliptic Painlev\'e equation is reviewed in ${\mathbb P}^1\times {\mathbb P}^1$ formalism.
In Section~\ref{sect:lax} the Lax pair for the elliptic Painlev\'e
equation is formulated. Some properties of relevant polynomials
are prepared in Section~\ref{sect:keys}. Finally,
in Section~\ref{sect:comp}, the compatibility
condition of the Lax pair is analyzed and its equivalence to the
elliptic Painlev\'e equation is established (Theorem~\ref{thm:main}).
In Appendix~\ref{sect:p6}, the dif\/ferential case is discussed.

Before closing this introduction, let us look at an observation
which may be helpful to motivate our construction.
In this paper, we see the Lax equations like
\eqref{eq:qP6L1}, \eqref{eq:qP6L2} from two dif\/ferent viewpoints.
One is a standard way, where we consider the equations as
dif\/ference equations for unknown function $y(z)$, and variables
$(f,g)$ are regarded as parameters. The other is unusual viewpoint,
where we consider these equations as equations of algebraic
curves in variables $(f,g) \in {\mathbb P}^1\times{\mathbb P}^1$, and
$y(z), y(qz), y(z/q)$ or $T^{-1}(y(z))$ are regarded as parameters.
In the second point of view, we have
\begin{proposition}
The equation \eqref{eq:qP6L1} is uniquely characterized as a curve
of bi-degree $(3,2)$ in ${\mathbb P}^1\times{\mathbb P}^1$ passing through the $12$ points:
\begin{gather}
(0,b_1/q), \quad (0,b_2/q), \quad (\infty,b_3), \quad (\infty, b_4),\quad
(a_1,0), \quad (a_2,0), \quad (a_3,\infty), \quad (a_4,\infty), \nonumber\\
(z,\infty),\quad  \left(\frac{z}{q},0\right),\quad
\left(z,\frac{a_1a_2}{q}\frac{y(z)}{y(qz)}\right),\quad
\left(\frac{z}{q},\frac{a_1a_2}{q}\frac{y(z/q)}{y(z)}\right). \label{eq:qP612pt}
\end{gather}
Similarly, the equation \eqref{eq:qP6L2} is also characterized as
a curve of bi-degree $(1,1)$ passing through $3$~points:
\begin{gather*}
(\infty,\infty),\quad  \left(z,\frac{a_1a_2}{q}\frac{y(z)}{y(qz)}\right),\quad
\left(z-\frac{a_1a_2}{z}\frac{y(z)}{T^{-1}y(z)},0\right).
\end{gather*}
\end{proposition}

In Sakai's theory, Painlev\'e equations are characterized by
the 9 points conf\/igurations in ${\mathbb P}^2$ or equivalently
by the 8 points conf\/igurations in ${\mathbb P}^1\times {\mathbb P}^1$.
We note that the f\/irst 8 points in~\eqref{eq:qP612pt} are nothing
but the conf\/iguration which characterize $q$-$P_{\rm VI}$.
This kind of relations between the Lax equations and the point
conf\/igurations have been observed also in other dif\/ference or
dif\/ferential cases~\cite{Y2} (See Appendix~\ref{sect:p6} for
$P_{\rm VI}$ case). Hence, it is naturally expected that the
Lax equations for the elliptic Painlev\'e equation will also
be determined by suitable conditions as plane algebraic curves.
This is what we will show in this paper.

\section{The elliptic Painlev\'e equation}\label{sect:ell-p}

Let $P_1, \dots, P_8$ be points on ${\mathbb P}^1 \times {\mathbb P}^1$.
We assume that the conf\/iguration of the
points $P_1, \dots, P_8$ is generic, namely the curve $C_0$ of bi-degree
(2,2) passing through the eight points is unique and it is a smooth
elliptic curve. We denote the
equation of the curve $C_0: \varphi_{22}(f,g)=0$, where $(f,g)$
is an inhomogeneous coordinate of ${\mathbb P}^1 \times {\mathbb P}^1$.
Let $X$ be the rational surface obtained from ${\mathbb P}^1 \times {\mathbb P}^1$
by blowing up the eight points $P_1, \dots, P_8$.
Its Picard lattice ${\rm Pic}(X)$ is given by
\begin{gather*}
{\rm Pic}(X)={\mathbb Z} H_1\oplus {\mathbb Z} H_2\oplus {\mathbb Z} E_1 \oplus \cdots \oplus {\mathbb Z} E_8,
\end{gather*}
where $H_i$ ($i=1,2$) is the class of lines corresponding to $i$-th
component of ${\mathbb P}^1 \times {\mathbb P}^1$ and  $E_j$ ($j=1,\ldots,8$) is the
exceptional divisors.
The nontrivial intersection pairings for these basis are given by
\begin{gather*}
(H_1,H_2)=(H_2,H_1)=1, \qquad (E_j,E_j)=-1.
\end{gather*}
Note that the surface $X$ is birational equivalent with the 9 points
blown-up of ${\mathbb P}^2$.

In the most generic situation, the group of Cremona transformations on
the surface $X$ is the af\/f\/ine Weyl group of type $E^{(1)}_8$ and its
translation part ${\mathbb Z}^8$ gives the elliptic Painlev\'e equations~\mbox{\cite{Sakai,Murata}}.
A choice of $E^{(1)}_8$ simple roots $\alpha_0, \ldots, \alpha_8$
in ${\rm Pic}(X)$ is
$\alpha_0=E_1-E_2$, $\alpha_1=H_1-H_2$, $\alpha_2=H_2-E_1-E_2$, $\alpha_i=E_{i-1}-E_i$ $(i=3,\ldots,8).$
The null root $\delta$ (=$-K_X$ : the anti-canonical divisor of $X$) is
\begin{gather}
\delta=2H_1+2H_2-E_1-E_2-\cdots-E_8.\label{eq:delPic}
\end{gather}
The action of the translation $T_{\alpha}$ on ${\rm Pic}(X)$ is
given by the Kac's formula
\begin{gather*}
T_{\alpha}(\beta)=\beta+(\delta,\beta)\alpha
-\left((\delta,\beta)\frac{(\alpha,\alpha)}{2}+(\alpha,\beta)\right)\delta.
\end{gather*}
For instance, for the translation $T=T_{E_i-E_j}$ along the direction
$E_i-E_j$ ($1\leq i\neq j \leq 8$), we have
\begin{gather}
T(H_i)=H_i+2(E_i-E_j)+2\delta, \qquad i=1,2,\nonumber \\
T(E_j)=E_i, \nonumber \\
T(E_i)=E_i+(E_i-E_j)+2 \delta,\nonumber \\
T(E_k)=E_k+(E_i-E_j)+\delta, \qquad k \neq i,j.\label{eq:TPic}
\end{gather}
Similar to the case of the 9 points blown-up of ${\mathbb P}^2$ \cite{KMNOY},
the above type of translations $T_{E_i-E_j}$ admit
simple geometric description as follows.

$(i)$ Points $P_1, \ldots, P_8$ are transformed as
\begin{gather}
T(P_k)=P_k, \qquad k \neq i,j,\nonumber \\
P_1+\cdots+P_{i-1}+T(P_i)+P_{i+1}+\cdots+P_8=0,\nonumber \\
T(P_i)+T(P_j)=P_i+P_j,\label{eq:Tpara}
\end{gather}
with respect to the addition on the elliptic curve $C_0$ passing through
$P_1, \ldots, P_8$.

$(ii)$
The transformation of the Painlev\'e dependent variable $P=(f,g)$ can be
found as follows. Let $C$ be the elliptic curve passing through
$P_1, \ldots ,P_{i-1}$, $P_{i+1}, \ldots, P_8$ and $P$.
It is easy to see that $T(P_i)$ lies on $C$. Def\/ine $T(P)$ by
\begin{gather}
T(P_i)+T(P)=P_j+P,\label{eq:TP}
\end{gather}
with respect to the addition on $C$.

In this paper, we employ the rules $(i)$ and $(ii)$ as the def\/inition
of the elliptic Painlev\'e equation. Then the relations~\eqref{eq:TPic}
are consequence of them.

It is convenient to introduce a Jacobian parametrization of the point
$P_u=(f_u,g_u)$ on $C_0$ in such a way that (1)~$P_u+P_v=P_{u+v}$, and
(2)~Let $C_{mn}$ be a curve of bi-degree $(m,n)$ and let
$P_{x_i}$ ($i=1,\ldots,2mn$) be the intersections $C_{mn}\cap C_0$, then
\begin{gather}
m h_1+n h_2-x_1-\cdots-x_{2mn}=0, \qquad ({\rm mod.} \ {\rm period}),
\label{eq:Abel}
\end{gather}
where $h_1$, $h_2$ are constant parameters.

We put\footnote{This $\delta \in {\mathbb C}$ is dif\/ferent from $\delta \in {\rm Pic}(X)$ in
equation~\eqref{eq:delPic}.} $\delta=2 h_1+2 h_2-u_1-\cdots-u_8$ where $u_i$ is the
parameter corresponding to the point $P_i=P_{u_i}$.
Note that $f_u=f_{h_1-u}$ and $g_u=g_{h_2-u}$.
An example of such parametrization is
\begin{gather*}
f_u=\frac{[u+a][u-h_1-a]}{[u+b][u-h_1-b]}, \qquad
g_u=\frac{[u+c][u-h_2-c]}{[u+d][u-h_2-d]},
\end{gather*}
where $[u]$ is an odd theta function and $a$, $b$, $c$, $d$ are constants.
An expression of the elliptic Painlev\'e equation on ${\mathbb P}^1\times {\mathbb P}^1$
using a parametrization in terms of the Weierstrass $\wp$ function
was given by Murata~\cite{Murata}.

In this paper, we will consider the case $T=T_{E_2-E_1}$ as an example,
and we use the notation:
\begin{gather*}
\dot{x}=T_{E_2-E_1}(x),
\end{gather*}
for any variables~$x$.
From equation~\eqref{eq:Tpara}, we have
\begin{gather*}
\dot{u_k}=u_k, \qquad k \neq 1,2,\qquad
\dot{u_1}=u_1-\delta,\qquad
\dot{u_2}=u_2+\delta.
\end{gather*}

In our construction, various polynomials and curves in
${\mathbb P}^1\times {\mathbb P}^1$ are def\/ined through their degree and vanishing
conditions. Let us introduce a notation to describe them.
\begin{definition}
Let $\Phi_{mn}(p_{1}^{m_1}p_{2}^{m_2}\cdots)$ be a
linear space of polynomials in $(f,g)$ of bi-degree $(m,n)$
which vanish at point $p_i\in {\mathbb P}^1\times {\mathbb P}^1$ with multiplicity $m_i$.
\end{definition}
Common zeros of $F \in \Phi_{mn}(p_{1}^{m_1}p_{2}^{m_2}\cdots)$
are called the base points of the family.
Note that there may be some un-assigned base points
besides to the assigned ones $p_1, p_2, \dots$.

For convenience, we also use an extended notation such as
\begin{gather*}
\Phi_{mn}^{d}\left(p_{1}^{m_1}p_{2}^{m_2}
\cdots\, |\, {p'}_1^{n_1}{p'}_2^{n_2}\cdots\right).
\end{gather*}
Where $d$ and ${p'}_1^{n_1}{p'}_2^{n_2}\cdots$ indicate the additional
information: the dimension
\begin{gather*}
d={\rm dim}\Phi_{mn}\left(p_{1}^{m_1}p_{2}^{m_2}\cdots p_{k}^{m_k}\,|\,\cdots\right)
=(m+1)(n+1)-\sum_{i=1}^k \frac{m_i(m_i+1)}{2},
\end{gather*}
and the un-assigned base points ${p'}_i$ with multiplicity $n_i$.

\section{The Lax equations}\label{sect:lax}

In this section we def\/ine a pair of 2nd order linear
dif\/ference equations (the Lax pair for the elliptic Painlev\'e
equation).

We chose a generic point $P_z$ on a curve $C_0$. The variable $z$
plays the role of dependent variable of the Lax equations.
Unknown function of the Lax equations is denoted by $y=y(z)$.
For simplicity, we use the following notation:
\begin{gather}
\overline{F}(z)=F(z+\delta), \qquad
\underline{F}(z)=F(z-\delta). \nonumber
\end{gather}
Then our Lax pair takes the following form:
\begin{gather}
{\rm (L1)}\quad L_1=C_1 \underline{y}+C_2 y+C_3 \overline{y}=0, \nonumber \\
{\rm (L2)}\quad L_2=C_4 \underline{y}+C_5 y+C_6 \dot{y}=0. \label{eq:lax12}
\end{gather}
Here $\dot{y}=T(y)$, and the coef\/f\/icients $C_1,\ldots,C_6$ depend on
$P_1,\ldots,P_8, P_z$ and $P=(f,g)$.

The main idea of our construction is to consider the equations like
\eqref{eq:lax12} as equations of curves
in variables $(f,g) \in {\mathbb P}^1\times{\mathbb P}^1$.

The f\/irst Lax equation {\rm (L1)} is def\/ined as follows.
\begin{definition}\label{def:L1}
Let $Q_z$ and $Q_{\underline{z}}$ be points in ${\mathbb P}^1\times{\mathbb P}^1$
def\/ined in the inhomogeneous
coordinates $(f,g)$ as
\begin{gather}
Q_z=\{f=f_z\}\cap\{(g-g_z)y=(g-g_{h_1-z})\overline{y}\},\nonumber \\
Q_{\underline{z}}=\{f=f_{\underline{z}}\}\cap\{(g-g_{\underline{z}})\underline{y}=(g-g_{h_1-\underline{z}})y\}.
\label{eq:Qz}
\end{gather}
Note that these points depend on $\overline{y}$, $y$, $\underline{y}$ besides the dependence
on $P_1, \ldots, P_8$ and $z$.
Then the curve $L_1=0$ is
def\/ined by the following conditions:
\begin{gather*}
{\rm (L1a)} \quad \mbox{$L_1 \in \Phi_{32}^3(P_1\cdots P_8P_z\, |\, P_{\delta+h_1-z})$},\\
{\rm (L1b)} \quad \mbox{the curve $L_1=0$ passes through $Q_z$ and $Q_{\underline{z}}$}.
\end{gather*}
\end{definition}

\begin{lemma}\label{lem:1}
The conditions {\rm (L1a)}, {\rm (L1b)} determine the curve $L_1=0$
uniquely and it is of the form {\rm (L1)} in equation~\eqref{eq:lax12}.
\end{lemma}

\begin{proof}
Polynomial of bi-degree $(3,2)$ has 12
free parameters.
The condition (L1a) determines~9 of them and we have 3 parameter
(2 dimensional) family of curves
\begin{gather*}
c_1 G_1(f,g)+c_2 G_2(f,g)+c_3 G_3(f,g)=0
\end{gather*}
satisfying the condition (L1a).
The condition (L1b) adds two more linear equations on the
coef\/f\/icients $c_1$, $c_2$, $c_3$, hence the curve $L_1=0$ is unique
up to an irrelevant overall factor. To see the resulting equation is
linear in $\underline{y}$, $y$, $\overline{y}$, we take the following
basis of the above family:
\begin{gather*}
G_1=(f-f_z){\varphi}_{22}(f,g), \qquad
G_2={\varphi}_{32}(f,g), \qquad
G_3=(f-f_{\underline{z}}){\varphi}_{22}(f,g).
\end{gather*}
Where, ${\varphi}_{22}=0$ is the equation of the curve $C_0$ and
${\varphi}_{32}$ is a polynomial of bi-degree (3,2) which is tangent to
the lines $f=f_z$ and $f=f_{\underline{z}}$ at $P_z$ and $P_{h_1-\underline{z}}$
respectively.
Then we have
\begin{gather*}
G_1=0, \qquad G_2\propto (g-g_z)^2, \qquad G_3\propto (g-g_z)(g-g_{h_1-z}),
\qquad {\rm for} \quad f=f_z,\nonumber \\
G_1\propto (g-g_{\underline{z}})(g-g_{h_1-\underline{z}}), \qquad G_2
\propto (g-g_{h_1-\underline{z}})^2, \qquad G_3=0,
\qquad {\rm for} \quad f=f_{\underline{z}},
\end{gather*}
and hence, $c_1 \propto \underline{y}$, $c_2 \propto y$, $c_3 \propto \overline{y}$.
\end{proof}

The 2nd Lax equation {\rm (L2)} in \eqref{eq:lax12} is def\/ined in a similar way.
\begin{definition}
Let $Q_{u_1}$ be a point on ${\mathbb P}^1\times {\mathbb P}^1$ given in inhomogeneous
coordinate $(f,g)$ as
\begin{gather}
Q_{u_1}=\{f=f_{u_1}\}\cap\{(g-g_{u_1})y=(g-g_{h_1-u_1})\dot{y}\},\label{eq:Q1}
\end{gather}
which depends on the variables $y$, $\dot{y}$.
Then the curve $L_2=0$ is def\/ined as

\begin{enumerate}\itemsep=0pt

\item[]{\rm (L2a)} \quad $L_2\in \Phi_{32}^3(P_1P_3\cdots P_8 P_{z+u_2-u_1}
P_{h_1+\delta-z}\,|\,P_1)$,\footnote{Due to \eqref{eq:Abel}, the intersection $(L_2=0)\cap C_0$ at
$P_1$ is of multiplicity~2. This means the curve $L_2=0$ is tangent to
$C_0$ at $P_1$, but is not a node in general.}

\item[]{\rm (L2b)} \quad
\mbox{the curve $L_2=0$ passes through $Q_{\underline{z}}$ in equation~\eqref{eq:Qz}
and $Q_{u_1}$}.
\end{enumerate}
\end{definition}

The fact that the curve specif\/ied above is unique and is of the form {\rm (L2)}
in equation~\eqref{eq:lax12} can be proved in a similar way as
Lemma \ref{lem:1}. In this case, {\rm (L2)} takes the form
\begin{gather*}
c_1 (f-f_1){\varphi}_{22}\underline{y}+c_2 F_{32}(h_1-\underline{z}) y
+c_3  (f-f_{\underline{z}}){\varphi}_{22}\dot{y}=0,
\end{gather*}
where the curve $F_{32}(h_1-\underline{z})=0$ is tangent to the lines $f=f_1$ and
$f=f_{\underline{z}}$ at $P_1$ and $P_{h_1-\underline{z}}$ respectively.
Then the curve $F_{32}(h_1-\underline{z})=0$ is tangent both
$f=f_1$ and $C_0$ at $P_1$, i.e.\ it has a~node at $P_1$. Hence
$F_{32}(h_1-\underline{z}) \in \Phi_{32}^2(P_1^2P_3\cdots P_8P_{h_1-\underline{z}})$.

In what follows, this polynomial $F_{32}(z)$ (a polynomial in $(f,g)$ with
parameter $z$) plays important role.
Its def\/ining properties are
\begin{gather}
F_{32}(z) \in \Phi_{32}^2(P_1^2P_3\cdots P_8P_{z}),\nonumber \\
F_{32}(z)=0 {\rm \ is \ tangent \ to \ the \ line\ }
f=f_{z} {\rm \ at \ } P_{z}.\label{eq:F32def}
\end{gather}
Under these conditions, $F_{32}(z)$ is unique up to normalization.

\section{Some useful relations}\label{sect:keys}

In this section, we prepare several formulas satisf\/ied by
$f$, $g$, $\dot{f}$ and $\dot{g}$. Some results (Lemmas~\ref{lem:54f},~\ref{lem:F1F2} and~\ref{lem:FoverF}) will be used to analyze
the compatibility of the Lax equations in the next section.

\begin{lemma}
For generic $Q=(x,y)\in {\mathbb P}^1\times{\mathbb P}^1$, let
$F=F(f,g)$ be a polynomial such that
$F \in \Phi_{54}^1(P_1^4P_3^2\cdots P_8^2Q)$.
Then $F=0$ for  $\dot{f}=\dot{x}$.
\end{lemma}
\begin{proof}
From equation \eqref{eq:TPic}, the evolution
$\dot{P}=(\dot{f},\dot{g})$ of $P=(f,g)$ takes the following form
\begin{gather*}
\dot{f}=\frac{F_1(f,g)}{F_2(f,g)}, \qquad
\dot{g}=\frac{G_1(f,g)}{G_2(f,g)}, 
\end{gather*}
where
$F_1, F_2 \in \Phi_{54}^2(P_1^4P_3^2\cdots P_8^2)$.
Then the polynomial $F\in \Phi_{54}^1(P_1^4P_3^2\cdots P_8^2Q)$
is given by $F \propto F_1(P)F_2(Q)-F_2(P)F_1(Q)$.
Then we have
$F=0 \Leftrightarrow \dot{f}=F_1(P)/F_2(P)=F_1(Q)/F_2(Q)=\dot{x}$
for $F_2(P)\neq 0$ and $F_2(Q)\neq 0$.
\end{proof}

From equation \eqref{eq:TP} we have
\begin{gather}
\dot{P}_z=P_{z+u_1-u_2-\delta}.\label{eq:evPz}
\end{gather}
Then, putting $Q=P_{z+\delta-u_1+u_2}$ (i.e. $\dot{Q}=P_z$)
in the above Lemma, we have
\begin{lemma}\label{lem:54f}
Let ${\varphi}_{54}(z) \in \Phi_{54}^1(P_1^4P_3^2\cdots
P_8^2P_{z+\delta-u_1+u_2}\, |\, P_{h_1+\delta-z-u_1+u_2})$.
Then ${\varphi}_{54}(z)=0$ for $\dot{f}=f_{z}$.
\end{lemma}

This lemma gives a characterization of $\dot{f}$ which
will be used in the next section (Lemma \ref{lem:L4}).
For the later use, we should also prepare a characterization
of $\dot{g}$ using some properties of the polynomial $F_{32}(z)$.
To do this, let us introduce an involution $r$ on ${\mathbb P}^1\times{\mathbb P}^1$:
\begin{gather}
r: \ (f, g) \mapsto \big(f, \tilde{g}(f,g)\big), \label{eq:r-def}
\end{gather}
def\/ined as follows.
For generic $Q=(x,y)$, let
$F(f,g) \in \Phi_{2,2}^1(P_1\dot{P}_2P_{3}\cdots P_8Q)$.
The equation $F(x,g)=0$ have two solutions,
one is trivial $g=y$, and the other solution $g=\tilde{g}(x,y)$
gives the desired birational transformation.
The action of the involution $r$ on the ${\rm Pic}(X)$ is given by
\begin{gather}
r(H_1)=H_1, \qquad
r(H_2)=4H_1+H_2-E_1-\cdots-E_8, \qquad
r(E_i)=H_1-E_i. \label{eq:rPic}
\end{gather}
Hence, $\tilde{g}(f,g)$ is a fractional linear transformation of $g$
with coef\/f\/icients depending on $f$.
Specialized to generic point on the curve $C_0$, we have
\begin{gather*}
r(P_{z})=P_{h_1-z}.
\end{gather*}

The basic property of the transformations $r$ and $T$ is
\begin{gather}\label{eq:rTPic}
r T(\lambda)=\lambda, \qquad \lambda=3H_1+2H_2-2E_1-E_3-\cdots- E_8,
\end{gather}
which follows from equations \eqref{eq:TPic} and \eqref{eq:rPic}.
More precisely, we have the following
\begin{lemma}
Let $\{F_1, F_2, F_3\}$ be a basis of polynomials
$\Phi_{3,2}^{3}(P_1^2P_3\cdots P_8)$, then the equation
$(P$ is given and $P'$ is unknown$)$
\begin{gather}
(F_1(P):F_2(P):F_3(P))=(F_1(P'):F_2(P'):F_3(P'))
\label{eq:F123}
\end{gather}
has unique unassigned solution $P'=rT(P)$.
\end{lemma}
\begin{proof} The equation \eqref{eq:F123} is equivalent to
$F_i(P)F_3(P')=F_3(P)F_i(P')$ ($i=1,2$), which are of bi-degree (3,2)
and have 12 solutions. 11 of them are assigned ones $P_1$ (multiplicity
$2^2=4$), $P_3, \ldots, P_8$ and trivial one $P'=P$, and hence there
exist one unassigned solution, which is given by $rT(P)$ by the
above formula \eqref{eq:rTPic}.
\end{proof}

The following is a special case of the above lemma.
\begin{lemma}\label{lem:F1F2}
Let $\{F_1, F_2\}$ be a basis of polynomials
$\Phi_{3,2}^{2}(P_1^2P_3\cdots P_8P_z)$, then we have
\begin{gather*}
T\left(\frac{F_1}{F_2}\right)=r\left(\frac{F_1}{F_2}\right), \qquad  \forall\; z.
\end{gather*}
\end{lemma}

In the remaining part of this section, we will study a special
polynomial ${\mathcal F}$. Its property (Lemma~\ref{lem:FoverF}) will play
crucial role in the next section.
As a polynomial in $(f,g)$, ${\mathcal F}$ is def\/ined by the conditions:
\begin{gather}
{\mathcal F}\in \Phi_{32}^2\left(P_1^2P_3\cdots P_8Q\right),\qquad
\frac{\partial {\mathcal F}}{\partial g}\Big{|}_{P=Q}=0. \label{eq:Fprop}
\end{gather}
Then ${\mathcal F}$ is unique up to normalization factor.
Note that the specialization ${\mathcal F}|_{Q=P_z}$ satisfy the def\/ining property of
$F_{32}(z)$ in equation~\eqref{eq:F32def}.
The normalization of ${\mathcal F}$ may depend on $Q$. We f\/ix it so that ${\mathcal F}$ is
a polynomial in $Q=(x,y)$ of minimal degree. Then we have
\begin{lemma}
As a polynomial in $Q=(x,y)$, ${\mathcal F}$ has bi-degree $(5,2)$ and
has zeros at $P_1$ $($double point$)$, $P_3, \dots, P_8$, $P$. Moreover it
satisfy the following properties:
\begin{gather}
\frac{\partial {\mathcal F}}{\partial y}\Big{|}_{Q=P_i}=0, \qquad  i=3,\ldots,8.
\nonumber
\end{gather}
\end{lemma}
\begin{proof}
Consider the following $12\times 12$ determinant:
\begin{gather}
D=m_{P_1}\wedge\frac{\partial m_{P_1}}{\partial f}\wedge
\frac{\partial m_{P_1}}{\partial g}\wedge m_{P_3}\wedge \cdots \wedge
m_{P_8}\wedge m_{P}\wedge m_{Q}\wedge \frac{\partial m_{Q}}{\partial y},
\label{eq:detD}
\end{gather}
where
\begin{gather*}
m_{(f,g)}=\left\{\left(1,f,f^2,f^3\right), \left(1,f,f^2,f^3\right)g, \left(1,f,f^2,f^3\right)g^2\right\}
\in {\mathbb C}^{12}
\nonumber
\end{gather*}
is a vector of monomials of bi-degree $(3,2)$.
As a polynomial in $(f,g)$, it is easy to see that
this determinant $D$ has the desired property \eqref{eq:Fprop} as ${\mathcal F}$.
As a polynomial in $Q=(x,y)$, the bi-degree of $D$ is apparently $(6,4)$.
The degree in variable $y$ is actually 2, since the $y$ dependent part
$m_{Q}\wedge \frac{\partial m_{Q}}{\partial y}$ in the determinant
can be reduced to
\begin{gather*}
\left\{\left(1,x,x^2,x^3\right), \left(1,x,x^2,x^3\right)\frac{y}{2}, (0,0,0,0)\right\}
\wedge\left\{(0,0,0,0), \left(1,x,x^2,x^3\right), \left(1,x,x^2,x^3\right)2y\right\}.
\end{gather*}
Moreover, the determinant $D$ is factorized by $(x-f_1)$ where $P_1=(f_1,g_1)$.
This follows form the relation
\begin{gather*}
2 m_{P_1}+(y-g_1)\frac{\partial m_{P_1}}{\partial g_1}-2 m_{Q}+
(y-g_1)\frac{\partial m_{Q}}{\partial y}=0 \qquad {\rm at} \quad x=f_1.
\end{gather*}
Hence, one can take ${\mathcal F}=D/(x-f_1)$ which is of degree $(5,2)$ in variables
$(x,y)$. Its desired vanishing conditions are easily checked from the
structure of the determinant \eqref{eq:detD}.
Since we have 17 vanishing conditions, the degree $(5,2)$ is minimal.
\end{proof}

\begin{lemma}\label{lem:Gprop}
For the determinant $D$ in \eqref{eq:detD}, we have
\begin{gather*}
D=(g-g_1)^2(x-f_1)^2 G \qquad {\rm at} \quad f=f_1. 
\end{gather*}
Where $G$ is independent of $P=(f,g)$ and is a polynomial in
$Q=(x,y)$ of degree $(4,2)$. It satisfy
\begin{gather*}
G=\frac{\partial G}{\partial y}=0 \qquad {\rm at}
\quad Q=P_1, P_3, \dots, P_8.
\end{gather*}
\end{lemma}

\begin{proof}
It is enough to show that
\begin{gather*}
D=\frac{\partial D}{\partial x}=0 \qquad {\rm at} \quad f=x=f_1.
\end{gather*}
To see this, let $M_i$ be the $i$-th vector in determinant
$D$ in equation~\eqref{eq:detD}.
Then for $f=x=f_1$ we have the following linear relations:
\begin{gather*}
(g-y)(g+y-2g_1)M_1+(g-y)(g-g_1)(y-g_1)M_3+(y-g_1)^2M_{10}
=(g-g_1)^2M_{11},  \\
2(g_1-y)M_1+(g-g_1)(g+g_1-2 y)M_3+2(y-g_1)M_{10}=(g-g_1)^2M_{12}.
\end{gather*}
Hence, $M_{11}\wedge M_{12}$ and $\frac{\partial}{\partial x}
(M_{11}\wedge M_{12})$
vanishes when multiplied with $M_1\wedge M_3\wedge M_{10}$.
\end{proof}

\begin{lemma}\label{lem:ABC1}
Let $G$ be the polynomial in the above Lemma {\rm \ref{lem:Gprop}} and
let $A=A(x)$, $B=B(x)$, $C=C(x)$ be the coefficient of the fractional
linear transformation:
\begin{gather}
\tilde{y}=\tilde{y}(x,y)=-\frac{A+By}{B+Cy}, \label{eq:fracABC}
\end{gather}
where $(x,\tilde{y})$ is the image of $(x,y)$ under the involution~$r$~\eqref{eq:r-def}.
Then we have $G=A+2 By+Cy^2$ up to a normalization factor.
\end{lemma}

\begin{proof}
Let $\phi_{22}(f,g){=}\phi_{22}(f,g;x,y)$ be a polynomial of $P{=}(f,g)$
belonging to $\Phi_{22}^1(P_1P_3{\cdots} P_8 Q)$ with $Q=(x,y)$.
By the def\/inition of the involution $r$, we have
\begin{gather*}
\phi_{22}(x,\tilde{y};x,y)=(y-\tilde{y})(A+B(y+\tilde{y})+C y\tilde{y}).
\end{gather*}
On the other hand, the polynomial $\phi_{22}(f,g;x,y)$ can be
represented as the following $9\times 9$ determinant:
\begin{gather*}
\phi_{22}(f,g;x,y)=m'_{P_1}\wedge m'_{P_3}\wedge
\cdots\wedge m'_{P_8}\wedge m'_{Q}\wedge m'_{P},
\end{gather*}
where $m'_{(f,g)} \in {\mathbb C}^9$ is a vector of monomials
$f^ig^j$ ($0\leq i,j \leq 2$).
Then we have
\begin{gather*}
A+2B y+C y^2
=
\mathop{\lim}_{\tilde{y}\rightarrow y}
\frac{\phi_{22}(x,\tilde{y};x,y)}{y-\tilde{y}}
=m'_{P_1}\wedge m'_{P_3}\wedge \cdots\wedge m'_{P_8}\wedge
m'_{Q}\wedge \frac{\partial m'_{Q}}{\partial y}.
\end{gather*}
The last determinant is of degree $(4,2)$ in $(x,y)$ and satisfy
the vanishing properties for $G$ in Lemma \ref{lem:Gprop}.
\end{proof}

\begin{lemma}\label{lem:ABC2}
For the polynomial $G=A+2By+Cy^2$ and
$\tilde{y}=\tilde{y}(x,y)$, $\tilde{g}=\tilde{g}(f,g)$, we have
\begin{gather*}
\frac{G(x,\tilde{y})}{G(x,y)}=\frac{AC-B^2}{(B+Cy)^2}
=\frac{\partial \tilde{y}}{\partial y}
=-\frac{(\tilde{g}-\tilde{y})(g-\tilde{y})}{(\tilde{g}-y)(g-y)}
\Big{|}_{f=x}.
\end{gather*}
\end{lemma}
\begin{proof}
All equalities follow
from direct computation by using the transformation \eqref{eq:fracABC}:
\begin{gather*}
 \tilde{y}(x,y)=-\frac{A+By}{B+Cy} \qquad \mbox{and}\qquad
  \tilde{g}(f,g)\Big{|}_{f=x}=-\frac{A+B g}{B+C g} .\tag*{\qed}
\end{gather*}\renewcommand{\qed}{}
\end{proof}

\begin{lemma}\label{lem:FoverF}
The following relation holds:
\begin{gather*}
\frac{{\mathcal F}(f,g:x,\tilde{y})}{{\mathcal F}(f,g;x,y)}\Big{|}_{f=f_1}
=-\frac{(\tilde{g}-\tilde{y})(g-\tilde{y})}{(\tilde{g}-y)(g-y)}
\Big{|}_{f=x},
\end{gather*}
and both sides of this equation are actually independent of $g$.
\end{lemma}
\begin{proof} This is a corollary of the
Lemmas \ref{lem:Gprop}, \ref{lem:ABC1} and \ref{lem:ABC2}.
\end{proof}

\section{The compatibility}\label{sect:comp}

The compatibility of the Lax pair (L1), (L2) in
equation \eqref{eq:lax12} is analyzed through the following
four steps (Fig.~\ref{fig:lax}).
\begin{enumerate}\itemsep=-1pt
\item Eliminating $\underline{y}$ from {\rm (L1)} and {\rm (
L2)} $\rightarrow$
equation {\rm (L3)} between $y$, $\overline{y}$, $\dot{y}$.
\item Eliminating $\underline{y}$ from {\rm (L2)} and \underline{\rm (L3)} $\rightarrow$
equation {\rm (L4)} between $y$, $\dot{y}$, $\underline{\dot{y}}$.
\item Eliminating $\overline{y}$ from $\overline{\rm (L2)}$ and {\rm (L3)} $\rightarrow$
equation {\rm (L5)} between $y$, $\dot{y}$, $\overline{\dot{y}}$.
\item Eliminating $y$ from {\rm (L4)} and ${\rm (L5)}$  $\rightarrow$
equation {\rm (L6)} between $\dot{y}$, $\underline{\dot{y}}$, $\dot{\overline{y}}$.
\end{enumerate}
Then the compatibility means the equivalence
{\rm (L6)} $\Leftrightarrow$ $T_{E_2-E_1}({\rm L1})$
which is the main result of this paper (Theorem \ref{thm:main}).

\begin{figure}[t]
\centering
\setlength{\unitlength}{1mm}
\begin{picture}(60,30)(0,5)
\put(10,10){$\underline{y}$}
\put(30,10){$y$}
\put(50,10){$\overline{y}$}
\put(10,30){$\dot{\underline{y}}$}
\put(30,30){$\dot{y}$}
\put(50,30){$\dot{\overline{y}}$}
\put(30,4){$L_1$}
\put(20,17){$L_2$}
\put(40,17){$L_3$}
\put(14,10){\line(1,0){15}}
\put(14,13){\line(1,0){15}}
\put(29,13){\line(0,1){15}}
\put(33,10){\line(1,0){15}}
\put(33,13){\line(1,0){15}}
\put(33,13){\line(0,1){15}}
\end{picture}
\hskip10mm
\begin{picture}(60,35)(0,5)
\put(10,10){$\underline{y}$}
\put(30,10){$y$}
\put(50,10){$\overline{y}$}
\put(10,30){$\dot{\underline{y}}$}
\put(30,30){$\dot{y}$}
\put(50,30){$\dot{\overline{y}}$}
\put(30,35){$L_6$}
\put(20,22){$L_4$}
\put(40,22){$L_5$}
\put(14,30){\line(1,0){15}}
\put(14,27){\line(1,0){15}}
\put(29,27){\line(0,-1){15}}
\put(33,30){\line(1,0){15}}
\put(33,27){\line(1,0){15}}
\put(33,27){\line(0,-1){15}}
\end{picture}
\caption{Lax equations.}
\label{fig:lax}
\end{figure}

We will track down these equations step by step.
The resulting properties are summarized 
as follows.
$$\begin{array}{|cc|l|l@{\,}|l@{\,}|}
\hline
{\rm equation}&{\rm term}&{\rm coef\/f\/icient}&{\rm divisors}
&{\rm additional}\ {\rm zeros}\\
\hline
   &\underline{y}&(f-f_{z}){\varphi}_{22}&(H_1)+(\delta)&P_{z},P_{h_1-z}\\
L_1&   y &         {\varphi}_{32}&(H_1+\delta) &P_{z},P_{h_1+\delta-z}\\
   &\overline{y}&(f-f_{\underline{z}}){\varphi}_{22}&(H_1)+(\delta)&P_{z-\delta},P_{h_1+\delta-z}\\
\hline
   &\underline{y}&(f-f_1){\varphi}_{22}  &(H_1-E_1)+(\delta)& -\\
L_2&   y &F_{32}(h_1-\underline{z})&(H_1+\delta-E_1+E_2)&P_{z-u_1+u_2},P_{h_1+\delta-z}\\
   &\dot{y}&(f-f_{\underline{z}}){\varphi}_{22}&(H_1)+(\delta)&P_{z-\delta},P_{h_1+\delta-z}\\
\hline
   &\dot{y}&(f-f_{z}){\varphi}_{22}&(H_1)+(\delta)&P_{z},P_{h_1-z}\\
L_3&   y &F_{32}(z)&(H_1+\delta-E_1+E_2)&P_{z},P_{h_1+\delta-u_1+u_2-z}\\
   &\overline{y}&(f-f_1){\varphi}_{22}  &(H_1-E_1)+(\delta)& -\\
\hline
   & y &{\varphi}_{54}(\underline{z})&(H_1+2\delta-2E_1+2E_2)&P_{z-u_1+u_2},P_{h_1+2\delta-u_1+u_2-z}\\
L_4&\dot{y} &F_{32}(\underline{z}){\varphi}_{22}&(H_1+\delta-E_1+E_2)+(\delta)&P_{z-\delta},P_{h_1+2\delta-u_1+u_2-z}\\
   &\dot{\underline{y}}&(f-f_1)({\varphi}_{22})^2  &(H_1-E_1)+2(\delta)& -\\
\hline
   & y &{\varphi}_{54}(z)&(H_1+2\delta-2E_1+2E_2)&P_{z+\delta-u_1+u_2},P_{h_1+\delta-u_1+u_2-z}\\
L_5&\dot{y} &F_{32}(h_1-z){\varphi}_{22}&(H_1+\delta-E_1+E_2)+(\delta)&P_{z+\delta-u_1+u_2},P_{h_1-z}\\
   &\dot{\overline{y}}&(f-f_1)({\varphi}_{22})^2  &(H_1-E_1)+2(\delta)& -\\
\hline
   &\dot{\underline{y}}&{\varphi}_{54}(z){\varphi}_{22}&(H_1+2\delta-2E_1+2E_2)+(\delta)&P_{z+\delta-u_1+u_2},P_{h_1+\delta-u_1+u_2-z}\\
L_6&\dot{y}&{\varphi}_{76}&(H_1+3\delta-2E_1+2E_2)&P_{z+\delta-u_1+u_2},P_{h_1+2\delta-u_1+u_2-z}\\
   &\dot{\overline{y}}&{\varphi}_{54}(\underline{z}){\varphi}_{22}&(H_1+2\delta-2E_1+2E_2)+(\delta)&P_{z-u_1+u_2},P_{h_1+2\delta-u_1+u_2-z}\\
\hline
\end{array}$$

\noindent Step 1:
\begin{lemma}
The Lax equation $L_3=0$ is uniquely characterised by the
following properties:
\begin{gather*}
{\rm (L3a)} \quad \mbox{$L_3 \in \Phi_{32}^3(P_1P_3\cdots P_8P_z
P_{h_1+\delta-z-u_1+u_2}\, |\, P_1)$},\\
{\rm (L3b)} \quad \mbox{passing through $2$ more points: $Q_{u_1}$ in \eqref{eq:Q1}
and $Q_{z}$ in \eqref{eq:Qz}}.
\end{gather*}
\end{lemma}

\begin{proof}
The property (L3b) follows directly from the corresponding conditions
in (L1b) and (L2b). Let us consider the property (L3a).
We know that the Lax equations (L1), (L2) have the following form:
\begin{gather*}
{\rm (L1)}\quad
L_1=(f-f_{z}){\varphi}_{22}\underline{y}+F y
+\ast (f-f_{\underline{z}}){\varphi}_{22}\overline{y}=0,  \\
{\rm (L2)}\quad
L_2=(f-f_1){\varphi}_{22}\underline{y}+F' y
+\ast (f-f_{\underline{z}}){\varphi}_{22}\dot{y}=0.
\end{gather*}
Here $F$, $F'$ are some polynomials of degree $(3,2)$
and $*$ represent some constant independent of~$(f,g)$.
From the equation $(f-f_1) L_1-(f-f_{z}) L_2=0$, we have
three term relation bet\-ween~$y$,~$\overline{y}$,~$\dot{y}$.
This relation is apparently of degree $(4,2)$, however, it is
divisible by $f-f_{\underline{z}}$. Since if it is not so, then it follows
that $y=0$ for $f=f_{\underline{z}}$ and for any $g$, which contradict the
2nd conditions of (L1b), (L2b). Then the quotient should belong to
$\Phi_{32}^3(P_1^2P_3\cdots P_8P_z\,|\,P_{h_1+\delta-z-u_1+u_2})$
as desired.
Uniqueness follows by a simple dimensional argument as before.
\end{proof}

The coef\/f\/icients in (L2), (L3) are related as follows.
\begin{lemma}\label{lem:A2B2}
For the normalized equations
\begin{gather*}
{\rm (L2)} \quad \underline{y}-A_2(z)y+B_2(z)\dot{y}=0,  \\
{\rm (L3)} \quad \overline{y}-A_3(z)y+B_3(z)\dot{y}=0,
\end{gather*}
we have
\begin{gather*}
A_3(h_1+\delta-z)=A_2(z), \qquad B_3(h_1+\delta-z)=B_2(z).
 \end{gather*}
\end{lemma}
\begin{proof} This is because that the characterization
properties {\rm (L2)} and {\rm (L3)} are related by
$\underline{y}\leftrightarrow \overline{y}$ and $z\leftrightarrow h_1+\delta-z$.
\end{proof}

\noindent Step 2:
\begin{lemma}\label{lem:L4}
The Lax equation $L_4=0$ has the following characterizing properties:
\begin{gather}
{\rm (L4a)} \quad \mbox{$L_4 \in \Phi_{54}^3(P_1^3P_3^2\cdots P_8^2P_{z-u_1+u_2}
P_{h_1+2\delta-z-u_1+u_2}\,|\,P_1)$},\nonumber\\
{\rm (L4b)} \quad  \mbox{passing through $2$ more points: $Q_{u_1}$ in \eqref{eq:Q1}
and $\dot{Q}_{\underline{z}}$ defined by}\nonumber\\
\phantom{{\rm (L4b)} \quad}{}
\dot{Q}_{\underline{z}}=\{\dot{f}=f_{\underline{z}}\}\cap\{
(\dot{g}-g_{h_1-\underline{z}})\dot{y}=(\dot{g}-g_{\underline{z}})\underline{\dot{y}}\}.
\label{eq:dotQuz}
\end{gather}
\end{lemma}

\begin{proof}
Eliminating $\underline{y}$ from (L2) and \underline{(L3)},
one get three term relation between $y$, $\underline{\dot{y}}$, $\dot{y}$.
It is apparently of degree $(6,4)$ but divisible by $f-f_{\underline{z}}$.
It is easy to check that quotient $L_4$ belongs to
$\Phi_{54}^3(P_1^3P_3^2\cdots P_8^2P_{z-u_1+u_2}
P_{h_1+2\delta-z-u_1+u_2}\,|\,P_1)$.

The f\/irst condition in (L4b) is the direct consequence of (L2b) or (L3b).

We will show the second condition in (L4b).
Using the Lemma~\ref{lem:A2B2}, the (L4) equation can be written as
\begin{gather*}
K y+A_2(z')B_2(z)\dot{y}+B_2(z')\dot{\underline{y}}=0,
\end{gather*}
where $z'=h_1+2\delta-z$ (i.e.\ $\underline{z}+\underline{z'}=h_1$)
and the coef\/f\/icient of $y$ is $K=1-A_2(z)A_2(z')$.
By tracing the zeros, we see that the numerator of $K$ is proportional to
\begin{gather*}
{\varphi}_{54} (\underline{z}) \in \Phi_{54}^1\left(P_1^4P_3^2\cdots P_8^2
P_{z-u_1+u_2}P_{h_1+2\delta-z-u_1+u_2}\right).
\end{gather*}
Hence, by Lemma \ref{lem:54f},
we have $K=0$ when $\dot{f}=f_{\underline{z}}$.
Thus, we have{\samepage
\begin{gather*}
\frac{\dot{\underline{y}}}{\dot{y}}=-b(z)A_2(z')
\qquad {\rm for} \quad \dot{f}=f_{\underline{z}}.
\end{gather*}
Here, we put
$b(z)=\frac{B_2(z)}{B_2(z')}$ which is independent of $(f,g)$.}

$A_2(z')|_{ \dot{f}=f_{\underline{z}}}$ is evaluated as follows.
By Lemma \ref{lem:F1F2}, we have
$A_2(z;f,g)=A_2(z;\dot{f},\tilde{\dot{g}})$. Hence,
by using the condition (L2b), we have
\begin{gather}
A_2(z)|_{\dot{f}=f_{\underline{z}}}=
\frac{\tilde{\dot{g}}-g_{\underline{z'}}}{\tilde{\dot{g}}-g_{\underline{z}}},
\qquad {\rm i.e.} \qquad
A_2(z')|_{\dot{f}=f_{\underline{z}}}=
\frac{\tilde{\dot{g}}-g_{\underline{z}}}{\tilde{\dot{g}}-g_{\underline{z'}}}.
\label{eq:A2-eval}
\end{gather}

Next, let us compute the factor $b(z)$. Since
\begin{gather*}
\mathop{\lim}_{f\rightarrow f_1}\frac{A_2(z)}{B_2(z)}=\frac{\dot{y}}{y}
=\frac{g-g_{u_1}}{g-g_{h_1-u_1}}
\end{gather*}
is independent of $z$, and
$A_2(z)=\frac{F_{32}(h_1-\underline{z})}{(f-f_1){\varphi}_{22}}$, we have
\begin{gather}
b(z)=\frac{B_2(z)}{B_2(z')}
=\frac{A_2(z)}{A_2(z')}\Big{|}_{f=f_1}
=\frac{F_{32}(\underline{z'})}{F_{32}(\underline{z})}\Big{|}_{f=f_1}.
\label{eq:b-eval1}
\end{gather}
Now, we apply the Lemma \ref{lem:FoverF} in case of $Q=(x,y)=P_{\underline{z}}$.
Then we have $x=f_{\underline{z}}=f_{\underline{z'}}$, $y=g_{\underline{z}}$,
$\tilde{y}=g_{\underline{z'}}$, ${\mathcal F}|_{Q=P_{\underline{z}}}\propto
F_{32}(\underline{z})$, and hence
\begin{gather}
\frac{F_{32}(\underline{z'})}{F_{32}(\underline{z})}\Big{|}_{f=f_1}
=\frac{{\mathcal F}(f,g:x,\tilde{y})}{{\mathcal F}(f,g;x,y)}\Big{|}_{f=f_1, Q=P_{\underline{z}}}
=-\frac{(\tilde{\dot{g}}-g_{\underline{z'}})(\dot{g}-g_{\underline{z'}})}
{(\tilde{\dot{g}}-g_{\underline{z}})(\dot{g}-g_{\underline{z}})}\Big{|}_{\dot{f}=f_{\underline{z}}}.
\label{eq:b-eval2}
\end{gather}
Here, in the last equation, variables $(f,g)$ is replaced
by $(\dot{f},\dot{g})$ using the $g$ independence of the expression.
It follows from \eqref{eq:A2-eval}, \eqref{eq:b-eval1} and
\eqref{eq:b-eval2} that
\begin{gather*}
\frac{\dot{\underline{y}}}{\dot{y}}
=\frac{(\dot{g}-g_{\underline{z'}})}{(\dot{g}-g_{\underline{z}})} \qquad
{\rm at} \quad \dot{f}=f_{\underline{z}}.
\end{gather*}
This is the desired 2nd relation in (L4b).
\end{proof}

\noindent Step 3:
\begin{lemma}
The Lax equation $L_5=0$ has the following characterizing properties:
\begin{gather}
{\rm (L5a)} \quad
\mbox{$L_5 \in \Phi_{54}^3(P_1^3P_3^2\cdots P_8^2P_{z+\delta-u_1+u_2}
P_{h_1+\delta-z-u_1+u_2}\,|\,P_1)$},\nonumber\\
{\rm (L5b)} \quad
\mbox{passing through $2$ more points: $Q_{u_1}$ in \eqref{eq:Q1} and
$\dot{Q}_{z}$ defined by}\nonumber\\
\phantom{{\rm (L5b)} \quad}{}
\dot{Q}_{z}=\{\dot{f}=f_{z}\}\cap\{
(\dot{g}-g_{h_1-z})\dot{\overline{y}}=(\dot{g}-g_{z})\dot{y}\}.
\label{eq:dotQz}
\end{gather}
\end{lemma}
The proof is omitted since it is almost the same as Step 2.

\medskip

\noindent Step 4:
\begin{lemma}\label{lem:L6}
The Lax equation $L_6=0$ has the following characterizing properties:
\begin{gather*}
{\rm (L6a)} \quad
\mbox{$L_6 \in \Phi_{76}^3(P_1^5P_2P_3^3\cdots P_8^3P_{z+\delta-u_1+u_2}
\,|\,P_{h_1+2\delta-z-u_1+u_2})$},\\
{\rm (L6b)} \quad
\mbox{passing through $2$ more points: $\dot{Q}_{z}$ in \eqref{eq:dotQz}
and $\dot{Q}_{\underline{z}}$ in \eqref{eq:dotQuz}.}
\end{gather*}
\end{lemma}

\begin{proof}
Eliminating $y$ from equations
\begin{gather*}
{\rm (L4)}\quad
{\varphi}_{54}(\underline{z})y+\ast A_{32}(\underline{z}){\varphi}_{22} \dot{y}
+\ast (f-f_1)({\varphi}_{22})^2 \underline{\dot{y}}=0, \\
{\rm (L5)}\quad
{\varphi}_{54}(z)y+\ast A_{32}(h_1-z){\varphi}_{22} \dot{y}
+\ast (f-f_1)({\varphi}_{22})^2 \overline{\dot{y}}=0,
\end{gather*}
we have ${\varphi}_{54}(z)L_4-{\varphi}_{54}(\underline{z})L_5=0$.
Which is apparently of degree $(10,8)$ but divisible by
\mbox{$(f-f_1){\varphi}_{22}$},
hence we have the equation of degree $(7,6)$. The vanishing
conditions follows from that of $L_4$ and $L_5$.
\end{proof}

Now we have the main result of this paper:
\begin{theorem}[The compatibility]\label{thm:main}
The equation {\rm (L6)} is equivalent with equation {\rm (L1)}
evolved by the translation $T_{E_2-E_1}$:
\begin{gather*}
u_i \mapsto \dot{u_i}, \qquad y \mapsto \dot{y},
\qquad (f,g) \mapsto (\dot{f}, \dot{g}).
\end{gather*}
Namely, the Lax pair {\rm (L1)}, {\rm (L2)} is compatible if and only if the
variables $(f,g)$ solve
the elliptic Painlev\'e equation for $T=T_{E_2-E_1}$.
\end{theorem}

\begin{proof}
We have obtained the characterization properties of (L6), hence
our task is to compare it with that for $T({\rm L1})$.

(1) From equation \eqref{eq:TPic}, we have
\begin{gather*}
T{\rm (L1)} \in T(H_1+\delta)=H_1+3\delta-2E_1+2E_2.
\end{gather*}

(2) Since (L1) has extra zeros at $P=P_z$ and $P=P_{h_1+\delta-z}$,
$T({\rm L1})$ has zeros at $\dot{P}=P_z$ and $\dot{P}=P_{h_1+\delta-z}$.
From the equation~\eqref{eq:evPz}, these extra zeros of $T{\rm (L1)}$ are
at $P=P_{z-u_1+u_2+\delta}$ and $P=P_{h_1+2\delta-z-u_1+u_2}$
in terms of original variable $P=(f,g)$.

From these two conditions, we see that $T{\rm (L1)}$ satisfy the
condition (L6a) in Lemma~\ref{lem:L6}. The condition~(L6b)
is exactly the condition~(L1b) transformed by~$T$.
\end{proof}

\appendix
\section[Discussion on the differential case]{Discussion on the dif\/ferential case}\label{sect:p6}

In this appendix, we will discuss the dif\/ferential case, taking
the sixth Painlev\'e equation $P_{\rm VI}$ as an example.
The $P_{\rm VI}$ equation has a Hamiltonian form
\begin{gather}
\frac{dq}{dt}=\frac{\partial H}{\partial p},\qquad
\frac{dp}{dt}=-\frac{\partial H}{\partial q}, \label{eq:p6}
\end{gather}
with Hamiltonian
\begin{gather}
H=\frac{1}{t(t-1)}\Big[q(q-1)(q-t)p^2+\big\{(a_1+2a_2)(q-1)q
+a_3(t-1)q+a_4t(q-1)\big\}p \nonumber\\
\phantom{H=}{} +a_2(a_1+a_2)(q-t)\Big].
\label{eq:p6ham}
\end{gather}
The equation \eqref{eq:p6} describes the iso-monodromy deformation
of the Fuchsian dif\/ferential equation
on ${\mathbb P}^1\setminus \{0,1,t,\infty\}$:
\begin{gather}
\frac{\partial^2 y}{\partial z^2}
+\left(\frac{1-a_4}{z}+\frac{1-a_3}{z-1}
+\frac{1-a_0}{z-t}-\frac{1}{z-q}\right)
\frac{\partial y}{\partial z}\nonumber \\
\phantom{\frac{\partial^2 y}{\partial z^2}}{}  +\left\{\frac{a_2(a_1+a_2)}{z(z-1)}
-\frac{t(t-1)H}{z(z-1)(z-t)}
+\frac{q(q-1)p}{z(z-1)(z-q)}\right\}y=0,
\label{eq:sl1-p6}
\end{gather}
($a_0+a_1+2a_2+a_3+a_4=1$) deformed by
\begin{gather}
\frac{\partial y}{\partial t}
+\frac{z(z-1)(q-t)}{t(t-1)(q-z)}\frac{\partial y}{\partial z}
+\frac{zp(q-1)(q-t)}{t(t-1)(z-q)}y=0.
\label{eq:sl2-p6}
\end{gather}
These equations \eqref{eq:sl1-p6}, \eqref{eq:sl2-p6} can be viewed
as a Lax pair for the $P_{\rm VI}$ equation.
To see the geometric meaning of these Lax equations,
let us introduce homogeneous coordinates $(X:Y:Z) \in {\mathbb P}^2$ by
\begin{gather}
q=\frac{Z}{Z-X}, \qquad p=\frac{Y(Z-X)}{XZ}.
\nonumber
\end{gather}
Then we have
\begin{proposition}
The equation \eqref{eq:sl1-p6} can be uniquely characterized as
an algebraic curve $F(X{,}Y{,}Z)\!$ $=0$
of degree~$4$ in ${\mathbb P}^2$, satisfying the following vanishing conditions:
\begin{gather*}
F(0,0,1)=F(1,-a_2,1)=F(1,0,0)=F(0,a_3,1)
=F(1,-a_1-a_2,1)=F(1,a_4,0)=0,
\\
F\left((t-1){\varepsilon},1,t{\varepsilon}-a_0t {\varepsilon}^2\right)=O\big({\varepsilon}^3\big),
\\
F\left((z-1){\varepsilon},1,z{\varepsilon}+z {\varepsilon}^2\right)=O\big({\varepsilon}^4\big),
\\
F\left(\frac{1}{z},\frac{1}{y}\frac{\partial y}{\partial z},
\frac{1}{z-1}\right)\Big{|}_{z \mapsto z+{\varepsilon}}=O\big({\varepsilon}^2\big).
\end{gather*}
Similarly the second Lax equation \eqref{eq:sl2-p6} has also a
similar characterization
as an algebraic curve $R(X,Y,Z)=0$ of degree $2$ with the following
conditions:
\begin{gather*}
R(0,1,0)=R(1,0,0)=0, \\
R\left((t-1){\varepsilon},1,t{\varepsilon}-\frac{t^2(t-z)}{z}\frac{1}{y}
\frac{\partial y}{\partial t}{\varepsilon}^2\right)=
O\big({\varepsilon}^3\big),  \\
R\left(\frac{1}{z},\frac{1}{y}\frac{\partial y}{\partial z},
\frac{1}{z}\right)=0.
\end{gather*}
\end{proposition}
This geometric characterization of the Lax equations for $P_{\rm VI}$
may be considered as a degenerate case of our construction.
The above result bear resemblance to the geometric characterization
of the Hamiltonian $H$ \eqref{eq:p6ham} as a cubic pencil~\cite{pen}.

Finally, let us give a comment on the apparent singularity and
the non-logarithmic property.
The Hamiltonian $H$ in equa\-tion~\eqref{eq:p6ham} is
usually f\/ixed by non-logarithmic condition for equa\-tion~\eqref{eq:sl1-p6}
at the apparent singularity $z=q$.
Namely, though the dif\/ferential equa\-tion~\eqref{eq:sl1-p6} has apparent singularity at $z=q$ with
exponents $0$ and $2$,
the solutions are actually holomorphic there.
In the dif\/ference Lax equation (L1) def\/ined in Def\/inition~\ref{def:L1},
the factor $(f-f_z)$ or $(f-f_{\underline{z}}$) in its coef\/f\/icients looks like
an ``apparent singularity''.
Since the non-logarithmicity is an essential property of the
dif\/ferential equation~\eqref{eq:sl1-p6},
it will be interesting if one can f\/ind
the corresponding notion in dif\/ference cases.

\subsection*{Acknowledgements}

The idea of this work came from the study of the Pad\'e
approximation method to the Painlev\'e equations \cite{Y}, and it was partially
presented at the Workshop
``Elliptic Integrable Systems, Isomonodromy Problems, and Hypergeometric Functions'' \cite{Y2}.
The author would like to thank the organisers and participants
for their interest.
He also thank to Professors K.~Kajiwara,
T.~Masuda, M.~Noumi, Y.~Ohta, H.~Sakai, M-H.~Saito and S.~Tsujimoto for discussions.
The author would like to thank the referees for their valuable comments and
suggestions.
This work is supported by Grants-in-Aid for Scientif\/ic No.17340047.

\pdfbookmark[1]{References}{ref}
\LastPageEnding

\end{document}